\documentclass[11pt]{article}
\usepackage{amsmath}
\usepackage{amssymb}
\def\R{\mathbf{R}}
\def\N{\mathbb{N}}

\newtheorem{theo}{\hspace*{\parindent}Theorem}
\newtheorem{lemma}{\hspace*{\parindent}Lemma}
\newtheorem{corol}{\hspace*{\parindent}Corollary}

\newcounter{theremark}
\newcommand{\rem}{\par\refstepcounter{theremark}\textbf{Remark \arabic{theremark}.} }
\title{Log-concavity for series in reciprocal gamma functions and applications}
\author{S.I.\,Kalmykov and D.B.\,Karp\footnote{Far Eastern Federal University and Institute of Applied Mathematics FEBRAS, Vladivostok,
Russia.  E-mails:\,\emph{sergeykalmykov@inbox.ru} (SK) and
\emph{dimkrp@gmail.com} (DK)}}
\date{}
\begin{document}
\maketitle

\begin{center}
\parbox{12cm}{
\small\textbf{Abstract.} Euler's  gamma function $\Gamma(x)$ is
logarithmically convex on $(0,\infty)$. Additivity of logarithmic
convexity implies that the function $x\to\sum{f_k\Gamma(x+k)}$ is
also log-convex (assuming convergence) if the coefficients are
non-negative. In this paper we investigate the series
$\sum{f_k\Gamma(x+k)^{-1}}$, where each term is clearly
log-concave. Log-concavity is not preserved by addition, so that
non-negativity of the coefficients is now insufficient to draw any
conclusions about the sum.  We demonstrate that the sum is
log-concave if $kf_k^2\geq(k+1)f_{k-1}f_{k+1}$ and is discrete
Wright log-concave if $f_k^2\geq{f_{k-1}f_{k+1}}$.  We conjecture
that the latter condition is in fact sufficient for the
log-concavity of the sum. We exemplify our general theorems by
deriving known and new inequalities for the modified Bessel,
Kummer and generalized hypergeometric functions and their
parameter derivatives.}
\end{center}

\bigskip


\bigskip


\bigskip

\paragraph{1. Introduction.}  A positive continuous
function $f$ defined on a real interval $I\subseteq\R$ is said to
be \textbf{Jensen log-concave} if
\begin{equation}\label{eq:log-conc}
f(\mu+h)^2\geq f(\mu)f(\mu+2h)
\end{equation}
for all $h>0$ and all $\mu$ such that $[\mu,\mu+2h]\subseteq{I}$.
If inequality (\ref{eq:log-conc}) is strict the function $f$ is
called strictly Jensen log-concave. If the sign of the inequality
is reversed one talks about Jensen log-convexity (or strict Jensen
log-convexity). For continuous functions Jensen log-concavity is
equivalent to log-concavity, i.e. concavity of $\log(f)$ (but is
weaker in general). It is also equivalent to the seemingly
stronger inequality
\begin{equation}\label{eq:wright-log-conc}
\phi_{h,s}(\mu)=f(\mu+h)f(\mu+s)-f(\mu)f(\mu+h+s)\geq 0 ~\text{for
all}~~h,s>0,
\end{equation}
which expresses the fact that $\mu\to f(\mu+h)/f(\mu)$ is
non-increasing for any $h>0$. We tacitly assume here and below
that all arguments lie in $I$. A function satisfying
(\ref{eq:wright-log-conc}) is called \textbf{Wright log-concave}
\cite[Chapter~I.4]{MPF}. For comparisons of these notions and
their higher order analogues see also the recent paper \cite{NRW}.

One is also frequently encountered with the situation when
(\ref{eq:log-conc}) or (\ref{eq:wright-log-conc}) only holds for
integer values of $h$.   We will express this fact by saying that
$f$ is \textbf{discrete log-concave} or \textbf{discrete Wright
log-concave}, respectively. In this case, however, we only have
the implication
(\ref{eq:wright-log-conc})$\Rightarrow$(\ref{eq:log-conc}), while
the reverse implication is not true even for continuous functions.
We note that $f$ is discrete Wright log-concave if and only if
\begin{equation}\label{eq:disc-wright}
\phi_{1,s}(x)=f(\mu+1)f(\mu+s)-f(\mu)f(\mu+s+1)\geq{0} ~\text{for
all}~~s>0.
\end{equation}
Indeed, (\ref{eq:disc-wright}) says that $f(\mu+1)/f(\mu)$ is
non-increasing, so that $f(\mu+2)/f(\mu+1)$ is again
non-increasing implying that their product $f(\mu+2)/f(\mu)$ is
non-increasing, i.e. satisfies (\ref{eq:wright-log-conc}) with
$h=2$.  In a similar fashion (\ref{eq:wright-log-conc}) holds for
all integer values of $h$ which is discrete Wright log-concavity
by definition. Discrete Jensen log-concavity and log-convexity are
also frequently referred to as ''Tur\'{a}n type inequalities''
following the classical result of Paul Tur\'{a}n for Legendre
polynomials \cite{Turan}.

If $f:\N_0\to{\R_+}$ is a sequence, then discrete log-concavity
reduces to the inequality $f_k^2\geq{f_{k-1}f_{k+1}}$, $k\in\N$.
We additionally require that the sequence $\{f_k\}_{k=0}^{\infty}$
is non-trivial and has no internal zeros:
$f_{N}=0~\Rightarrow~f_{N+i}=0$ for all $i\in\N_{0}$.  Such
sequences are also known as $PF_2$ (P\'{o}lya frequency sub two)
or doubly positive \cite{Karlin}.

Clearly, if $f$ is (Jensen or Wright) log-concave then $1/f$ is
(Jensen or Wright) log-convex. Notwithstanding the simplicity of
this relation, several important properties of log-concavity and
log-convexity differ. Particularly important is that log-convexity
is additive while log-concavity is not. Further, log-convexity is
a stronger property than convexity whereas log-concavity is weaker
than concavity.

In \cite{KS} the second author and Sergei Sitnik initiated the
investigation of the following problem: under what conditions on
non-negative sequence $\{f_k\}$ and the numbers $a_i$, $b_j$ the
function
\begin{equation}\label{eq:problem1}
\mu\to
f(\mu;x):=\sum\limits_{k=0}^{\infty}f_k\frac{\prod_{i=1}^{n}\Gamma(a_i+\mu+\varepsilon_ik)}
{\prod_{j=1}^{m}\Gamma(b_j+\mu+\varepsilon_{n+j}k)}x^k
\end{equation}
is (discrete) log-concave or log-convex?  Here $\Gamma$ is Euler's
gamma function and $\varepsilon_l$ can be $1$ or $0$. In \cite{KS}
the authors treated the cases of (\ref{eq:problem1}) with $n=1$,
$m=0$, $\varepsilon_1=1$; $n=m=1$, $\varepsilon_1=1$,
$\varepsilon_2=0$; and $n=m=1$, $\varepsilon_1=0$,
$\varepsilon_2=1$.  Of course, the log-convexity cases are nearly
trivial but the results of \cite{KS} go beyond log-convexity by
showing that the function $\phi_{h,s}(\mu;x)$ on the left-hand
side of (\ref{eq:wright-log-conc}) has non-positive Taylor
coefficients in powers of $x$, so that $x\to-\phi_{h,s}(\mu;x)$ is
absolutely monotonic. According to Hardy, Littlewood and P\'{o}lya
theorem \cite[Proposition~2.3.3]{NP} this also implies that this
function is multiplicatively concave:
$\phi_{h,s}(\mu;x^{\lambda}y^{1-\lambda})\geq{\phi_{h,s}(\mu;x)^{\lambda}\phi_{h,s}(\mu;y)^{1-\lambda}}$
for $\lambda\in[0,1]$.

In this paper we treat the case $n=0$, $m=1$, $\varepsilon_1=1$.
We get slightly different results depending on conditions imposed
on the sequence $\{f_{k}\}$.  If $f_{k}$ is log-concave without
internal zeros (i.e. doubly positive) we prove discrete Wright
log-concavity of the series (\ref{eq:problem1}). We conjecture
that the true log-concavity holds but we were unable to
demonstrate this.  If $\{f_kk!\}$ is doubly positive we show that
(\ref{eq:problem1}) is log-concave. We do so by establishing
non-negativity of the Taylor coefficients of $\phi_{h,s}(\mu;x)$
in powers of $x$ (either for all or only for integer $h>0$).
Again, by Hardy, Littlewood and P\'{o}lya theorem this implies
that $x\to\phi_{h,s}(\mu;x)$ is multiplicatively convex.

The paper is organized as follows: in section~2 we collect several
lemmas repeatedly used in the proofs; section~3 comprises two
theorems constituting the main content of the paper; in section~4
we give applications to Bessel, Kummer and generalized
hypergeometric functions and relate them to several previously
known results.

\paragraph{2. Preliminaries.}
We will need several lemmas which we prove in this section. First,
we formulate an elementary inequality we will repeatedly use
below.
\begin{lemma}\label{lm:uvrs}
Suppose $u,v,r,s>0$, $u=\max(u,v,r,s)$ and $uv>rs$. Then
$u+v>r+s$.
\end{lemma}
\textbf{Proof.} Indeed, dividing by $u$ we can rewrite the
required inequality as $r'+s'<1+v'$, where $r'=r/u\in(0,1)$,
$s'=s/u\in(0,1)$, $v'=v/u\in(0,1)$.  Since $r's'<v'$ by $rs<uv$,
the required inequality follows from the elementary inequality
$r'+s'<1+r's'$. ~~$\square$

Lemma~\ref{lm:uvrs} is a particular case of a much more general
result on logarithmic majorization - see \cite[2.A.b]{MOA}.

In the next lemma we say that a sequence has no more than one
change of sign if the pattern is
$(--\cdots--00\cdots00++\cdots++)$, where zeros and minus signs
may be omitted.
\begin{lemma}\label{lm:sum} Suppose
$\{f_k\}_{k=0}^{n}$ has no internal zeros and
$f_k^2\geq{f_{k-1}f_{k+1}}$, $k=1,2,\ldots,n-1$. If the real
sequence $A_0,A_1,\ldots,A_{[n/2]}$ satisfying $A_{[n/2]}>0$ and
$\sum\limits_{0\leq{k}\leq{n/2}}\!\!\!\!A_k\geq{0}$ has no more
than one change of sign, then
\begin{equation}\label{eq:keysum}
\sum\limits_{0\leq{k}\leq{n/2}}f_{k}f_{n-k}A_k\geq{0}.
\end{equation}
Equality is only attained if $f_k=f_0\alpha^k$, $\alpha>0$, and
$\sum\limits_{0\leq{k}\leq{n/2}}\!\!\!\!A_k=0$.
\end{lemma}
\textbf{Proof.} Suppose $f_k>0$, $k=s,\ldots,p$, $s\geq{0}$,
$p\leq{n}$. Log-concavity of $\{f_k\}_{k=0}^{n}$ clearly implies
that $\{f_{k}/f_{k-1}\}_{k=s+1}^{p}$ is decreasing, so that for
$s+1\leq{k}\leq{n-k+1}\leq{p+1}$
$$
\frac{f_k}{f_{k-1}}\geq\frac{f_{n-k+1}}{f_{n-k}}~\Leftrightarrow~f_{k}f_{n-k}\geq
f_{k-1}f_{n-k+1}.
$$
Since $k\leq{n-k+1}$ is true for all $k=1,2,\ldots,[n/2]$, the
weights $f_{k}f_{n-k}$ assigned to negative $A_k$s in
(\ref{eq:keysum}) are smaller than those assigned to positive
$A_k$s leading to (\ref{eq:keysum}).  The equality statement is
obvious.~~~$\square$
\begin{lemma}\label{lm:gammasum}
Suppose $m$ is non-negative integer. The inequality
\begin{equation}\label{eq:gammasum}
\sum\limits_{k=0}^{m}\frac{1}{k!(m-k)!}
\left(\frac{1}{\Gamma(k+\mu+a)\Gamma(m-k+\mu+b)}-\frac{1}{\Gamma(m-k+\mu)\Gamma(k+\mu+a+b)}\right)\geq{0}
\end{equation}
holds for  all  $a,b\geq{0}$, $\mu\geq{-1}$ and such that
$\mu+a\geq{0}$, $\mu+b\geq{0}$ . Equality is only attained if
$ab=0$.
\end{lemma}
\textbf{Proof.}  Using the easily verifiable relations
$$
(c)_k=\frac{\Gamma(c+k)}{\Gamma(c)},~~(m-k)!=\frac{(-1)^km!}{(-m)_k}~~\text{and}~~(c)_{m-k}=\frac{(-1)^k(c)_m}{(1-c-m)_k}
$$
we obtain
\begin{multline*}
\sum\limits_{k=0}^{m}\frac{1}{k!(m-k)!\Gamma(k+\alpha)\Gamma(m-k+\beta)}=\frac{1}{m!\Gamma(\alpha)\Gamma(\beta)}
\sum\limits_{k=0}^{m}\frac{(-1)^k(-m)_k}{k!(\alpha)_k(\beta)_{m-k}}
\\
=\frac{1}{m!\Gamma(\alpha)\Gamma(\beta)(\beta)_m}\sum\limits_{k=0}^{m}\frac{(-m)_k(1-\beta-m)_k}{k!(\alpha)_k}
=\frac{(\alpha+\beta+m-1)_m}{m!\Gamma(\alpha)(\alpha)_m\Gamma(\beta)(\beta)_m}
\\
=\frac{\Gamma(\alpha+\beta+2m-1)}{\Gamma(\alpha+m)\Gamma(\beta+m)\Gamma(\alpha+\beta+m-1)m!},
\end{multline*}
where we have used the Chu-Vandermonde identity
\cite[Corollary~2.2.3]{AAR}
$$
\sum\limits_{k=0}^{m}\frac{(-m)_k(a)_k}{(c)_kk!}=\frac{(c-a)_m}{(c)_m}.
$$
 This leads to an explicit evaluation of the left hand side of
(\ref{eq:gammasum}) as
$$
\frac{\Gamma(2\mu+a+b+2m-1)}{\Gamma(2\mu+a+b+m-1)m!}\left(
\frac{1}{\Gamma(m+\mu+a)\Gamma(m+\mu+b)}-\frac{1}{\Gamma(m+\mu)\Gamma(m+\mu+a+b)}
\right)
$$
For $a,b,\mu>0$, the required inequality reduces to log-convexity
of $\Gamma(x)$ for $x>0$ \cite[Corollary~1.2.6]{AAR}.  If $ab=0$
we clearly get the equality. If $a,b>0$, $\mu=m=0$, the second
term in parentheses disappears and (\ref{eq:gammasum}) holds
strictly. If $m=0$, $-1\leq\mu<0$ and $\mu+a\geq{0}$,
$\mu+b\geq{0}$ then the  term outside the parentheses reduces to 1
while the second term inside parentheses is strictly negative
(since $\mu+a+b>0$), so that the sum is strictly positive.  If
$m\geq{1}$ then $\mu+m\geq{0}$ and we are back in the previous
situation.~~$\square$

\begin{lemma}\label{lm:rec-gammas}
Suppose $m\geq{0}$ is an integer. Then for all complex $\beta$ and
$\mu$
\begin{multline}\label{eq:rec-gamma-sum}
S_m(\mu,\beta):=\sum\limits_{k=0}^{m}\left\{\frac{1}{\Gamma(k+\mu+1)\Gamma(m-k+\mu+\beta)}
-\frac{1}{\Gamma(k+\mu)\Gamma(m-k+\mu+\beta+1)}\right\}
\\
=\frac{(\mu+\beta)_{m+1}-(\mu)_{m+1}}{\Gamma(\mu+m+1)\Gamma(\mu+\beta+m+1)}.
\end{multline}
\end{lemma}
\textbf{Proof.}  We have
\begin{multline*}
S_m(\mu,\beta)=\frac{1}{\Gamma(\mu+1)\Gamma(\mu+\beta)\Gamma(\mu)\Gamma(\mu+\beta+1)}
\sum\limits_{k=0}^{m}\left\{\frac{\Gamma(\mu)\Gamma(\mu+\beta+1)}{(\mu+1)_k(\mu+\beta)_{m-k}}
-\frac{\Gamma(\mu+1)\Gamma(\mu+\beta)}{(\mu)_{k}(\mu+\beta+1)_{m-k}}\right\}
\\
=\frac{1}{\Gamma(\mu+1)\Gamma(\mu+\beta+1)}
\sum\limits_{k=0}^{m}\left\{\frac{\mu+\beta}{(\mu+1)_k(\mu+\beta)_{m-k}}
-\frac{\mu}{(\mu)_{k}(\mu+\beta+1)_{m-k}}\right\}
\\
=\frac{1}{\Gamma(\mu+1)\Gamma(\mu+\beta+1)}
\sum\limits_{k=0}^{m}\frac{1}{(\mu)_k(\mu+\beta)_{m-k}}
\left\{\frac{\mu(\mu+\beta)}{\mu+k}
-\frac{\mu(\mu+\beta)}{\mu+\beta+m-k}\right\}
\\
=\frac{1}{\Gamma(\mu+1)\Gamma(\mu+\beta+1)}
\sum\limits_{k=0}^{m}\frac{\beta+m-2k}{(\mu+1)_k(\mu+\beta+1)_{m-k}}.
\end{multline*}
Writing
$$
u_k=\frac{1}{(\mu+1)_{k-1}(\mu+\beta+1)_{m-k}},~1\leq{k}\leq{m},~~
u_0=\frac{\mu}{(\mu+\beta+1)_{m}},~u_{m+1}=\frac{\mu+\beta}{(\mu+1)_{m}},
$$
we get a telescoping sum, since
$$
u_{k+1}-u_{k}=\frac{1}{(\mu+1)_{k}(\mu+\beta+1)_{m-k-1}}-\frac{1}{(\mu+1)_{k-1}(\mu+\beta+1)_{m-k}}
=\frac{\beta+m-2k}{(\mu+1)_k(\mu+\beta+1)_{m-k}}
$$
for $k=0,1,\ldots,m$, so that
$$
\sum\limits_{k=0}^{m}(u_{k+1}-u_{k})=u_{m+1}-u_{0}=\frac{\mu+\beta}{(\mu+1)_{m}}-\frac{\mu}{(\mu+\beta+1)_{m}}
$$
and
$$
\frac{1}{\Gamma(\mu+1)\Gamma(\mu+\beta+1)}
\sum\limits_{k=0}^{m}\frac{\beta+m-2k}{(\mu+1)_k(\mu+\beta+1)_{m-k}}=
\frac{(\mu+\beta)_{m+1}-(\mu)_{m+1}}{\Gamma(\mu+m+1)\Gamma(\mu+\beta+m+1)}.~\square
$$
The following is a straightforward consequence of formula
(\ref{eq:rec-gamma-sum}).
\begin{corol}\label{cr:rec-gammas}
If $\mu\geq{-1}$, $\beta\geq{0}$ and $\mu+\beta\geq{0}$ then
$S_{m}(\mu,\beta)\geq{0}$.  The inequality is strict unless
$\beta=0$.
\end{corol}

\paragraph{3. Main results.} In this section we prove two general
theorems for series in reciprocal gamma functions.  The power
series expansions in this section  are understood as formal, so
that no questions of convergence are discussed. In applications
the radii of convergence will usually be apparent. The results of
this section are exemplified by concrete special functions in the
subsequent section.
\begin{theo}\label{th:gammadenom}
Suppose $\{f_n\}_{n=0}^{\infty}$ is a non-trivial non-negative
log-concave sequence without internal zeros. Then the function
\begin{equation}\label{eq:func}
\mu\mapsto
f(\mu,x)=\sum\limits_{n=0}^{\infty}\frac{f_nx^n}{n!\Gamma(\mu+n)},
\end{equation}
is strictly log-concave on $(0,\infty)$ for each fixed $x\geq{0}$.
Moreover, the function
$$
\varphi_{a,b,\mu}(x):=f(a+\mu,x)f(b+\mu,x)-f(a+b+\mu,x)f(\mu,x)=\sum\limits_{m=0}^{\infty}\varphi_mx^m
$$
has positive power series coefficients $\varphi_m>0$ for
$\mu\geq{-1}$, $a,b>0$ and  $\mu+a\geq{0}$, $\mu+b\geq{0}$ so that
$\varphi_{a,b,\mu}(x)$ is absolutely monotonic and
multiplicatively convex on $(0,\infty)$.
\end{theo}

\noindent\textbf{Proof.} Cauchy product yields
$$
\varphi_m=\sum\limits_{k=0}^{m}\frac{f_kf_{m-k}}{k!(m-k)!}\left(\frac{1}{\Gamma(k+\mu+a)\Gamma(m-k+\mu+b)}-\frac{1}{\Gamma(m-k+\mu)\Gamma(k+\mu+a+b)}\right).
$$
 Further, by Gauss summation (the first term is combined with the last,
 the second with the  second last etc.) we can write $\varphi_m$ in the form
\begin{equation}\label{eq:repr}
\varphi_m=\sum\limits_{k=0}^{[m/2]}\frac{f_kf_{m-k}}{k!(m-k)!}M_k(a,b,\mu),
\end{equation}
where for $k<m/2$
\begin{multline*}
M_k(a,b,\mu)=\underbrace{[\Gamma(k+\mu+a)\Gamma(m-k+\mu+b)]^{-1}}_{=u}
+\underbrace{[\Gamma(k+\mu+b)\Gamma(m-k+\mu+a)]^{-1}}_{=v}
\\[0pt]
-\underbrace{[\Gamma(m-k+\mu)\Gamma(k+\mu+a+b)]^{-1}}_{=r}
-\underbrace{[\Gamma(k+\mu)\Gamma(m-k+\mu+a+b)]^{-1}}_{=s}
\end{multline*}
and for $k=m/2$
$$
M_k=[\Gamma(m/2+\mu+a)\Gamma(m/2+\mu+b)]^{-1}-[\Gamma(m/2+\mu)\Gamma(m/2+\mu+a+b)]^{-1}.
$$
Under assumptions on $\mu,a,b$ made in the theorem
\begin{equation}\label{pos}
\sum\limits_{k=0}^{[m/2]}\frac{M_k(a,b,\mu)}{k!(m-k)!}\geq{0}
\end{equation}
according to Lemma~\ref{lm:gammasum}.  Write $M_k:=M_k(a,b,\mu)$
for brevity. We aim to show that the sequence
$\{M_k\}_{k=0}^{[m/2]}$ has no more than one change of sign with
$M_{[m/2]}>0$ in order to apply Lemma~\ref{lm:sum}. Due to
log-convexity of $\Gamma(x)$ the ratio $x\mapsto
\Gamma(x+\alpha)/\Gamma(x)$ is strictly increasing on $(0,\infty)$
when $\alpha>0$. This immediately implies $M_{[m/2]}>0$ and the
following inequalities
\begin{equation}\label{vu}
v\geq u \Leftrightarrow
\frac{\Gamma(m-k+\mu+b)}{\Gamma(m-k+\mu+a)}\geq
\frac{\Gamma(k+\mu+b)}{\Gamma(k+\mu+a)}~-~\text{true
for}~k\leq{m-k},
\end{equation}
\begin{equation}\label{us}
u > s \Leftrightarrow
\frac{\Gamma(m-k+\mu+a+b)}{\Gamma(m-k+\mu+b)}>
\frac{\Gamma(k+\mu+a)}{\Gamma(k+\mu)}~-~\text{true
for}~k\leq{m-k},
\end{equation}
\begin{equation}\label{vr}
v\geq r \Leftrightarrow \frac{\Gamma(k+\mu+a+b)}{\Gamma(k+\mu+b)}
\geq \frac{\Gamma(m-k+\mu+a)}{\Gamma(m-k+\mu)}~-~\text{true
for}~(m-b)/2\leq k,
\end{equation}
\begin{equation}\label{rv}
r\geq v \Leftrightarrow
\frac{\Gamma(m-k+\mu+a)}{\Gamma(m-k+\mu)}\geq
\frac{\Gamma(k+\mu+a+b)}{\Gamma(k+\mu+b)}~-~\text{true
for}~k\leq(m-b)/2.
\end{equation}
If $(m-b)/2\leq k<m/2$ then the sum of (\ref{us}) and (\ref{vr})
yields $M_k=u+v-r-s>0$. If $k<(m-b)/2$ then it follows from
(\ref{vu}), (\ref{us}) and (\ref{rv}) that $r>v>u>s$ (equality
cannot be attained in (\ref{vu}) and (\ref{rv}) under this
restriction on $k$). We will change notation to simplify writing:
$\alpha = \mu, \beta = b + \mu, \delta = a$. According the
hypothesis of the theorem we have $\beta>\alpha>0, \, \delta>0$.
We will show now that if $M_k<0 \Leftrightarrow u+v<r+s$ for some
$0<k<(m-b)/2$ then $M_{k-1}<0$ as well. Indeed, using
$z\Gamma(z)=\Gamma(z+1)$ we can write $M_{k-1}$ in the following
form
$$
M_{k-1}(\delta)=
\frac{k-1+\alpha+\delta}{m-k+\beta}u+\frac{k-1+\beta}{m-k+\alpha+\delta}v-\frac{k-1+\beta+\delta}{m-k+\alpha}r-\frac{k-1+\alpha}{m-k+\beta+\delta}s.
$$
Treating $u, v, r, s$ as constants we see that $M_{k-1}(0)<0$ by
forming a combination of $r>v$ and $r+s>v+u$ with positive
coefficients: $(C_1+C_2)r+C_2s>(C_1+C_2)v+C_2u$, where
$$
C_1+C_2=\frac{k-1+\beta}{m-k+\alpha}>C_2=\frac{k-1+\alpha}{m-k+\beta}.
$$
Further, differentiating with respect to $\delta$, we get
$$
M'_{k-1}(\delta)=
\frac{1}{m-k+\beta}u-\frac{k-1+\beta}{(m-k+\alpha+\delta)^2}v-\frac{1}{m-k+\alpha}r+\frac{k-1+\alpha}{(m-k+\beta+\delta)^2}s,
$$
which is obviously negative since $r>u$ (by (\ref{rv}) and
(\ref{vu})) and $v>s$ (by (\ref{vu}) and (\ref{us})). This shows
that $M_{k-1}< 0$ and hence $\{M_k\}_{k=0}^{[m/2]}$ has no more
than one change of sign. Applying Lemma~\ref{lm:sum} with
$A_k=M_k/(k!(m-k)!)$ we conclude that $\varphi_m>0$.
Multiplicative convexity follows by  Hardy, Littlewood and
P\'{o}lya theorem \cite[Proposition~2.3.3]{NP}. ~~$\square$

\rem If $\{f_n\}_{n=0}^{\infty}$ is log-convex then $\varphi_m$
can take both signs.

\begin{corol}\label{cr:compl-mon}
Assume the series in $(\ref{eq:func})$ converges for all
$x\geq{0}$. Then $\varphi_{a,b,\mu}(1/y)$ is completely monotonic
and log-convex on $[0,\infty)$, so that there exists a
non-negative measure $\tau$ supported on $[0,\infty)$ such that
\begin{equation}\label{eq:compl-mon}
\varphi_{a,b,\mu}(x)=\int\limits_{[0,\infty)}e^{-t/x}d\tau(t).
\end{equation}
\end{corol}
\textbf{Proof}. According to \cite[Theorem~3]{MS} a convergent
series of completely monotonic with non-negative coefficients is
again completely monotonic. This implies that
$y\to\varphi_{a,b,\mu}(1/y)$ is completely monotonic, so that  the
above integral representation follows by Bernstein's theorem
\cite[Theorem~1.4]{SSV}.  Log-convexity follows from complete
monotonicity according to \cite[Exersice 2.1(6)]{NP}.~~$\square$

In the next corollary  we adopt the convention
$\Gamma(-1)=-\infty$, $\Gamma(0)=+\infty$.
\begin{corol}\label{cr:f-twosided}
Under hypotheses and notation of Theorem~\ref{th:gammadenom}
\begin{equation}\label{eq:f-twosided}
\frac{\Gamma(a+\mu)\Gamma(b+\mu)}{\Gamma(\mu)\Gamma(a+b+\mu)}<\frac{f(\mu,x)f(a+b+\mu,x)}{f(a+\mu,x)f(b+\mu,x)}<1~\text{for
all}~x\geq{0}.
\end{equation}
If $\mu=0$ or $\mu=-1$ we additionally require that $x\ne0$
otherwise the left inequality becomes equality.
\end{corol}
\textbf{Proof.} The estimate from above is a restatement of
Theorem~\ref{th:gammadenom} since it is equivalent to
$\phi_{\beta,\mu}(x)>0$.

The estimate from below is obvious for $\mu=-1$ for we have zero
or negative number (if $a=1$ or $b=1$) on the left and a positive
number on the right for $x>0$. The remaining proof will be divided
into two cases (I) $\mu\geq{0}$; and (II)$-1<\mu<0$,
$\mu+a\geq{0}$, $\mu+b\geq{0}$ (recall that $a,b>0$ by hypotheses
of Theorem~\ref{th:gammadenom}).

In case (I) the left-hand inequality in (\ref{eq:f-twosided})
follows from strict log-convexity of $\mu\to\Gamma(\mu)f(\mu,x)$
which has been proved in \cite[Theorem~3]{KS} (where one has to
take account of the formula $(\mu)_k=\Gamma(\mu+k)/\Gamma(\mu)$).

In case (II) $\Gamma(\mu)<0$ and the left-hand inequality in
(\ref{eq:f-twosided}) reduces to
$$
\Gamma(a+\mu)f(a+\mu,x)\Gamma(b+\mu)f(b+\mu,x)>\Gamma(\mu)f(\mu,x)\Gamma(a+b+\mu)f(a+b+\mu,x).
$$
This inequality follows by observing that
$\Gamma(\mu)f(\mu,x)=\sum_{n=0}^{\infty}f_nx^n(n!(\mu)_n)^{-1}$
and
$$
\sum\limits_{k=0}^{m}\frac{f_{k}f_{m-k}}{k!(m-k)!}
\left\{\frac{1}{(a+\mu)_k(b+\mu)_{m-k}}-\frac{1}{(\mu)_k(a+b+\mu)_{m-k}}\right\}>0,
$$
since for $k=1,2,\ldots,m$ $(\mu)_k<0$ and for $k=0$
$(b+\mu)_{m}<(a+b+\mu)_{m}$.~~$\square$

\begin{corol}\label{cr:phi-below}
Under hypotheses and notation of Theorem~\ref{th:gammadenom} and
for all $x\geq{0}$
$$
f(a+\mu,x)f(b+\mu,x)-f(a+b+\mu,x)f(\mu,x)\geq
f_0^2\left[\frac{1}{\Gamma(\mu+a)\Gamma(\mu+b)}-\frac{1}{\Gamma(\mu)\Gamma(\mu+a+b)}\right]
$$
with equality only at $x=0$.
\end{corol}
\textbf{Proof.} Indeed, the claimed inequality is just
$\varphi_{a,b,\mu}(x)\geq\varphi_{a,b,\mu}(0)$ which is true by
Theorem~\ref{th:gammadenom}.~~$\square$

\rem Corollaries~\ref{cr:f-twosided} and \ref{cr:phi-below} imply
by elementary calculation the following bounds for the so called
``generalized Turanian''
$\Delta_{\varepsilon}(\mu,x):=f(\mu,x)^2-f(\mu+\varepsilon,x)f(\mu-\varepsilon,x)$:
\begin{equation}\label{eq:genTuranian}
A_\varepsilon(\mu)f_0^2\leq\Delta_{\varepsilon}(\mu,x)\leq
B_\varepsilon(\mu)f(\mu,x)^2,
\end{equation}
where $\mu-\varepsilon\geq{-1}$, $\mu\geq{0}$,  $x\geq{0}$ and
\begin{equation}\label{eq:AB}
A_\varepsilon(\mu)=\frac{\Gamma(\mu-\varepsilon)\Gamma(\mu+\varepsilon)-\Gamma(\mu)^2}
{\Gamma(\mu-\varepsilon)\Gamma(\mu+\varepsilon)\Gamma(\mu)^2},~~~~
B_\varepsilon(\mu)=\frac{\Gamma(\mu-\varepsilon)\Gamma(\mu+\varepsilon)-\Gamma(\mu)^2}
{\Gamma(\mu-\varepsilon)\Gamma(\mu+\varepsilon)}.
\end{equation}
In particular, if $\varepsilon=1$ the bounds
(\ref{eq:genTuranian}) simply to ($\mu\geq{0}$,  $x\geq{0}$)
\begin{equation}\label{eq:Turanian}
\frac{f_0^2}{\mu\Gamma(\mu)^2}\leq
f(\mu,x)^2-f(\mu+1,x)f(\mu-1,x)\leq
\frac{1}{\mu}f(\mu,x)^2,~~~x\geq{0},~\mu\geq{0}.
\end{equation}

Theorem~\ref{th:gammadenom} can be reformulated in terms of the
numbers $g_n:=f_n/n!$.  The hypotheses of the theorem require then
that these numbers satisfy
$$
g_n^2\geq\frac{n+1}{n} g_{n-1}g_{n+1}
$$
- a condition stronger then log-concavity.  If we weaken it to
log-concavity  we are only able to prove discrete Wright
log-concavity of $\mu\to\,f(\mu,x)$ in the next theorem. We
conjecture below that the adjective ``discrete'' is actually
redundant.
\begin{theo}\label{th:gammadenom1}
Suppose $\{g_n\}_{n=0}^{\infty}$ is a non-trivial non-negative
log-concave sequence without internal zeros. Then the function
\begin{equation}\label{eq:g-def}
\mu\to
g(\mu,x)=\sum\limits_{n=0}^{\infty}\frac{g_nx^n}{\Gamma(\mu+n)},
\end{equation}
is strictly discrete Wright log-concave on $(0,\infty)$ for each
fixed $x\geq{0}$. Moreover, the function
$$
\lambda_{\beta,\mu}(x):=g(\mu+1,x)g(\mu+\beta,x)-g(\mu,x)g(\mu+\beta+1,x)=\sum\limits_{m=0}^{\infty}\lambda_mx^m
$$
has positive power series coefficients $\lambda_m>0$ for each
$\mu\geq{-1}$ and $\beta>0$ such that $\mu+\beta\geq{0}$. This
implies that $x\to\lambda_{\beta,\mu}(x)$ is absolutely monotonic
and multiplicatively convex on $(0,\infty)$.
\end{theo}
\textbf{Proof.} Pursuing the same line of argument as in
Theorem~\ref{th:gammadenom} we have by the Cauchy product and the
Gauss summation:
$$
\lambda_m=\sum\limits_{k=0}^{[m/2]}g_kg_{m-k}M_k(1,\beta,\mu),
$$
where the numbers $M_k$ are defined in the proof of
Theorem~\ref{th:gammadenom}, below formula (\ref{eq:repr}). Under
assumptions on $\mu$ and $\beta$ made in the theorem we have
$$
\sum\limits_{k=0}^{[m/2]}M_k(1,\beta,\mu)=S_m(\mu,\beta)>0
$$
according to Corollary~\ref{cr:rec-gammas}. Further, it has been
shown in the course of the proof of Theorem~\ref{th:gammadenom}
that the sequence $\{M_k\}_{k=0}^{[m/2]}$ has no more than one
change of sign with $M_{[m/2]}>0$. Hence by Lemma~\ref{lm:sum} we
conclude that $\lambda_m>0$ implying discrete Wright log-concavity
of $\mu\to{g(\mu,x)}$ for each $x\geq{0}$ and absolutely
monotonicity of $x\to\lambda_{\beta,\mu}(x)$. Multiplicative
convexity follows by  Hardy, Littlewood and P\'{o}lya theorem
\cite[Proposition~2.3.3]{NP}. ~~$\square$

\begin{corol}\label{cr:g-compl-mon}
Assume the series in $(\ref{eq:g-def})$ converges for all
$x\geq{0}$. Then $\lambda_{\beta,\mu}(1/y)$ is completely
monotonic and log-convex on $[0,\infty)$, so that there exists a
non-negative measure $\tau$ supported on $[0,\infty)$ such that
\begin{equation}\label{eq:g-compl-mon}
\lambda_{\beta,\mu}(x)=\int\limits_{[0,\infty)}e^{-t/x}d\tau(t).
\end{equation}
\end{corol}

\begin{corol}\label{cr:g-twosided}
Under hypotheses and notation of Theorem~\ref{th:gammadenom1}
\begin{equation}\label{eq:g-twosided}
\frac{\mu}{\beta+\mu}<\frac{g(\mu,x)g(1+\beta+\mu,x)}{g(1+\mu,x)g(\beta+\mu,x)}<1~\text{for
all}~x\geq{0}.
\end{equation}
\end{corol}
\textbf{Proof.} The estimate from above is a restatement of
Theorem~\ref{th:gammadenom1} since it is equivalent to
$\lambda_{\beta,\mu}(x)>0$.

The estimate from below is obvious for $\mu=-1$ for we a negative
number on the left and a positive number on the right. The
remaining proof will be divided into two cases (I) $\mu\geq{0}$;
and (II)$-1<\mu<0$, $\mu+\beta\geq{0}$ (recall that $\beta>0$ by
hypotheses of Theorem~\ref{th:gammadenom1}).

In case (I) the left-hand inequality in (\ref{eq:g-twosided})
follows from strict log-convexity of $\mu\to\Gamma(\mu)g(\mu,x)$
which, in view of $(\mu)_k=\Gamma(\mu+k)/\Gamma(\mu)$, has been
proved in \cite[Theorem~3]{KS}.

In case (II)  the left-hand inequality in (\ref{eq:g-twosided})
can be rewritten as
$$
\Gamma(1+\mu)g(1+\mu,x)\Gamma(\beta+\mu)g(\beta+\mu,x)>\Gamma(\mu)g(\mu,x)\Gamma(1+\beta+\mu)g(1+b+\mu,x).
$$
This inequality follows by observing that
$\Gamma(\mu)g(\mu,x)=\sum_{n=0}^{\infty}g_n(\mu)_n^{-1}$ and
$$
\sum\limits_{k=0}^{m}g_{k}g_{m-k}
\left\{\frac{1}{(1+\mu)_k(\beta+\mu)_{m-k}}-\frac{1}{(\mu)_k(1+\beta+\mu)_{m-k}}\right\}>0,
$$
since  $(\mu)_k<0$ for $k=1,2,\ldots,m$  and
$(\beta+\mu)_{m}<(1+\beta+\mu)_{m}$ for $k=0$.~~$\square$
\begin{corol}\label{cr:lambda-below}
Under hypotheses and notation of Theorem~\ref{th:gammadenom1} and
for all $x\geq{0}$
$$
g(\mu+1,x)g(\mu+\beta,x)-g(\mu,x)g(\mu+\beta+1,x)\geq
\frac{g_0^2\beta}{\Gamma(\mu+1)\Gamma(\mu+\beta+1)}
$$
with equality only at $x=0$.
\end{corol}

\rem Since we have only proved discrete Wright log-concavity in
Theorem~\ref{th:gammadenom1} we cannot make any statements about
the ``generalized Turanian''
$g(\mu,x)^2-g(\mu+\varepsilon,x)g(\mu-\varepsilon,x)$. We can
assert, however, that the standard Turanian satisfies the
following bounds similar to those in (\ref{eq:Turanian})
\begin{equation}\label{eq:Turanian1}
\frac{g_0^2}{\mu\Gamma(\mu)^2}\leq
g(\mu,x)^2-g(\mu+1,x)g(\mu-1,x)\leq
\frac{1}{\mu}g(\mu,x)^2,~~~x\geq{0},~\mu\geq{0}.
\end{equation}

\textbf{Conjecture~1.} Under hypotheses of
Theorem~\ref{th:gammadenom1} the function $\mu\to{g(\mu,x)}$ is
log-concave on $(0,\infty)$ for each fixed $x\geq{0}$.  Moreover,
the function
$$
x\to g(\mu+\alpha,x)g(\mu+\beta,x)-g(\mu,x)g(\mu+\alpha+\beta,x)
$$
has positive power series coefficients for $\mu\geq{-1}$ and
$\mu+\alpha\geq{0}$, $\mu+\beta\geq{0}$, where $\alpha,\beta>0$.

The above conjecture is equivalent to the assertion that
\begin{equation*}
\sum\limits_{k=0}^{m}\left\{\frac{1}{\Gamma(k+\mu+\alpha)\Gamma(m-k+\mu+\beta)}
-\frac{1}{\Gamma(k+\mu)\Gamma(m-k+\mu+\alpha+\beta)}\right\}>0
\end{equation*}
which extends Corollary~\ref{cr:rec-gammas}.

\paragraph{4. Applications and relation to other work.}
We start with the well-studied case of the modified Bessel
function.  Even for this classical case we can add to the current
knowledge.

\medskip

\textbf{Example~1}. The modified Bessel function is defined by the
series \cite[formula~(4.12.2)]{AAR}
$$
I_{\nu}(u)=\sum\limits_{n\geq0}\frac{(u/2)^{2n+\nu}}{n!\Gamma(n+\nu+1)}.
$$
Hence, if we set $f_n=1$ $\forall n$, $x=(u/2)^{2}$ and
$\mu=\nu+1$ in Theorem~\ref{th:gammadenom} and use
$\partial^2_\nu\log(u/2)^{\nu}=0$ we immediately conclude that
$\nu\to{I_{\nu}(u)}$ is log-concave on $(-1,\infty)$ for each
fixed $u>0$. Moreover, for any $\nu\geq{-1}$ and
$\nu-\varepsilon\geq{-2}$ the ``generalized Turanian''
\begin{equation}\label{eq:I1}
u\to\Delta_{\varepsilon}(\nu,u):=(I_{\nu}(u))^2-I_{\nu+\varepsilon}(u)I_{\nu-\varepsilon}(u)
\end{equation}
has positive power series coefficients, is multiplicatively convex
and according to (\ref{eq:genTuranian}) satisfies
\begin{equation}\label{eq:I2}
(u/2)^{2\nu}A_\varepsilon(\nu+1)\leq\Delta_{\varepsilon}(\nu,u)\leq
B_\varepsilon(\nu+1)(I_{\nu}(u))^2,~~u\geq{0},
\end{equation}
where $A_\varepsilon$ and $B_\varepsilon$ are defined in
(\ref{eq:AB}).  In particular for $\varepsilon=1$ we get for
$\nu\geq{-1}$:
\begin{equation}\label{eq:I3}
\frac{(u/2)^{2\nu}}{(\nu+1)\Gamma(\nu+1)^2}\leq(I_{\nu}(u))^2-I_{\nu+1}(u)I_{\nu-1}(u)
\leq\frac{1}{\nu+1}(I_{\nu}(u))^2.
\end{equation}
All the more, the function
$(u/2)^{-2\nu}\Delta_{\varepsilon}(\nu,u)$ admits representation
(\ref{eq:compl-mon}). Proofs of various forms of log-concavity of
$I_{\nu}(u)$ have a long history. The discrete log-concavity,
$I_{\nu-1}(x)I_{\nu+1}(x)\leq{(I_{\nu}(x))^2}$, and the right-hand
side of (\ref{eq:I3}) for $\nu\geq{0}$ were probably first
demonstrated in 1951 by Thiruvenkatachar and Nanjundiah \cite{TN}.
In fact, our method here is an extension of their approach, so
that they could have proved the log-concavity of $\nu\to
I_{\nu}(x)$. The discrete log-concavity was rediscovered by Amos
\cite{Am} in 1974 and later by Joshi and Bissu \cite{JB} in 1991
with different proofs. Their paper also gives a proof of the
right-hand side of (\ref{eq:I3}) for $\nu\geq{0}$.  Finally, Lorch
in \cite{Lo} and later Baricz in \cite{Baricz10H} showed the
log-concavity of $\nu\to I_{\nu}(x)$ on $(-1,\infty)$ and
demonstrated the positivity of the function (\ref{eq:I1}) for
$\nu>-1/2$ and small $\varepsilon$. He also conjectured that the
positivity remains true for $\nu>-1$ and $\varepsilon\in(0,1]$.
Baricz \cite{Baricz10} demonstrated the Lorch's conjecture for
$\varepsilon=1$ and extended the right-hand side of (\ref{eq:I3})
to $\nu>-1$. Our results here not only confirm Lorch's conjecture
but also refine and strengthen it by proving (\ref{eq:I2}) and the
positivity of the power series coefficients of (\ref{eq:I1}).
Various extensions and a related results can also be found in
\cite{Baricz10,Baricz12,Segura}. We note that many proofs use
special properties of the modified Bessel functions, like
differential-recurrence relations, zeros etc.
Theorem~\ref{th:gammadenom} and its corollaries show that it is in
fact the structure of the power series that is responsible for the
bounds (\ref{eq:I2}) and (\ref{eq:I3}).

\medskip

\textbf{Example~2}.  In his 1993 preprint \cite{Sitnik} Sitnik,
among other things, proved the inequality
$$
R_{n}^2(x)>R_{n-1}(x)R_{n+1}(x), ~~x>0,~~n=1,2,\ldots,
$$
where
$$
R_{n}(x)=e^{x}-\sum\limits_{k=0}^{n}\frac{x^k}{k!}=\frac{x^{n+1}}{(n+1)!}{_1F_1}(1;n+2;x)
$$
is the exponential remainder.  We can generalize this function as
follows
$$
R_{\eta,\nu}(x)=\frac{x^{\nu+1}}{\Gamma(\nu+2)}{_1F_1}(\eta;\nu+2;x)
=x^{\nu+1}\sum\limits_{k=0}^{\infty}\frac{(\eta)_kx^k}{\Gamma(\nu+2+k)k!}.
$$
It is straightforward to check that the sequence $g_k=(\eta)_k/k!$
is log-concave iff $\eta\geq{1}$.  Then according to
Theorem~\ref{th:gammadenom1} the function $\nu\to R_{\eta,\nu}(x)$
is discrete Wright log-concave on $(-2,\infty)$ for each fixed
$\eta\geq{1}$, $x>0$ and
$$
x\to
R_{\eta,\nu+1}(x)R_{\eta,\nu+\beta}(x)-R_{\eta,\nu}(x)R_{\eta,\nu+\beta+1}(x)
$$
has positive power series coefficients for $\nu\geq-3$,
$\nu+\beta\geq-2$, $\beta>0$. Moreover,
$$
\frac{x^{2\nu+2}}{(\nu+2)\Gamma(\nu+2)^2}\leq
R_{\eta,\nu}(x)^2-R_{\eta,\nu+1}(x)R_{\eta,\nu-1}(x)\leq
\frac{1}{\nu+2}R_{\eta,\nu}(x)^2,~~~x\geq{0},~\nu\geq{-2}.
$$

\medskip

\textbf{Example~3}.  In addition to the results for the Kummer
function presented in Example~2 above we can derive bounds for its
logarithmic derivative.  The logarithmic derivatives of the Kummer
function plays an important role in some probabilistic
applications - see \cite{SK}.  Let us use abbreviated notation
$F(a;b;x)={_1F_1}(a;b;x)$.  The following contiguous relations are
easy to check (recall that $F'(a;b;x)=(a/b)F(a+1;b+1;x)$):
\begin{equation}\label{eq:con1}
aF(a;b;x)-aF(a+1;b;x)+xF'(a;b;x)=0,
\end{equation}
\begin{equation}\label{eq:con2}
abF(a+1;b;x)=b(a+x)F(a;b;x)-(b-a)xF(a;b+1;x),
\end{equation}
\begin{equation}\label{eq:con3}
b(b-1)(F(a;b-1;x)-F(a;b;x))-axF(a+1;b+1;x)=0.
\end{equation}
Dividing (\ref{eq:con2}) by $b$ and substituting $aF(a+1;b;x)$
into (\ref{eq:con1}) we get after simplification and dividing by
$\Gamma(b+1)$:
$$
\frac{1}{\Gamma(b+1)}F(a;b+1;x)=\frac{F(a;b;x)-F'(a;b;x)}{(b-a)\Gamma(b)}
$$
From (\ref{eq:con3}) we obtain:
$$
\frac{1}{\Gamma(b-1)}F(a;b-1;x)=\frac{1}{\Gamma(b)}(xF'(a;b;x)+(b-1)F(a;b;x))
$$
Thus we get the following expression for the Turanian:
\begin{multline*}
\frac{F(a;b;x)^2}{\Gamma(b)^2}-\frac{F(a;b+1;x)F(a;b-1;x)}{\Gamma(b+1)\Gamma(b-1)}
\\
=\frac{1}{\Gamma(b)^2(b-a)}\left\{-(a-1)F(a;b;x)^2+xF'(a;b;x)^2+(b-x-1)F(a;b;x)F'(a;b;x)\right\}
\end{multline*}
Hence, inequality (\ref{eq:Turanian}) becomes ($x>0$, $b>0$,
$a\geq{1}$):
$$
\frac{1}{b\Gamma(b)^2}<
\frac{-(a-1)F(a;b;x)^2+xF'(a;b;x)^2+(b-x-1)F(a;b;x)F'(a;b;x)}{\Gamma(b)^2(b-a)}
<\frac{1}{b\Gamma(b)^2}F(a;b;x)^2,
$$
which leads to ($F\equiv{F(a;b;x)}$, $F'\equiv{F'(a;b;x)}$):
$$
0<\frac{-(a-1)+x(F'/F)^2+(b-x-1)(F'/F)}{(b-a)}<\frac{1}{b}.
$$
Solving these quadratic inequalities we arrive at
$$
\frac{x+1-b+\sqrt{(x+1-b)^2+4x(a-1)}}{2x}<\frac{F'(a;b;x)}{F(a;b;x)}
<\frac{x+1-b+\sqrt{(x+1-b)^2+4xa(b-1)/b}}{2x}
$$
for $x>0$ and $b>a\geq{1}$.  The upper and lower bounds
interchange if $x>0$, $a\geq{1}$ and $0<b<a$:
$$
\frac{x+1-b+\sqrt{(x+1-b)^2+4xa(b-1)/b}}{2x}<\frac{F'(a;b;x)}{F(a;b;x)}
<\frac{x+1-b+\sqrt{(x+1-b)^2+4x(a-1)}}{2x}.
$$
These bounds are quite precise numerically especially when $a$ and
$b$ are close.

\medskip

\textbf{Example~4}.  The generalized hypergeometric function is
defined by the series
\begin{equation}\label{eq:pFqdefined}
{_{p}F_q}\left(\left.\!\!\begin{array}{c} a_1,a_2,\ldots,a_p\\
b_1,b_2,\ldots,b_q\end{array}\right|z\!\right):=\sum\limits_{n=0}^{\infty}\frac{(a_1)_n(a_2)_n\cdots(a_{p})_n}{(b_1)_n(b_2)_n\cdots(b_q)_nn!}z^n,
\end{equation}
where $(a)_0=1$, $(a)_n=a(a+1)\cdots(a+n-1)$, $n\geq{1}$, denotes
the rising factorial. The series (\ref{eq:pFqdefined}) converges
in the entire complex plane if $p\leq{q}$ and in the unit disk if
$p=q+1$.  In the latter case its sum can be extended analytically
to the whole complex plane cut along the ray $[1,\infty)$
\cite[Chapter~2]{AAR}. Applications of
Theorems~\ref{th:gammadenom} and \ref{th:gammadenom1} to
generalized hypergeometric function is largely based on the
following lemma.
\begin{lemma}\label{lm:HVVKS}
Denote by $e_k(x_1,\ldots,x_q)$  the $k$-th elementary symmetric
polynomial,
$$
e_0(x_1,\ldots,x_q)=1,~~~e_k(x_1,\ldots,x_q)=\!\!\!\!\!\!\!\!\sum\limits_{1\leq{j_1}<{j_2}\cdots<{j_k}\leq{q}}
\!\!\!\!\!\!\!\!x_{j_1}x_{j_2}\cdots{x_{j_k}},~~k\geq{1}.
$$
Suppose $q\geq{1}$ and $0\leq{r}\leq{q}$ are integers, $a_i>0$,
$i=1,\ldots,q-r$, $b_i>0$, $i=1,\ldots,q$, and
\begin{equation}\label{eq:symmetric-chain1}
\frac{e_q(b_1,\ldots,b_q)}{e_{q-r}(a_1,\ldots,a_{q-r})}\leq
\frac{e_{q-1}(b_1,\ldots,b_q)}{e_{q-r-1}(a_1,\ldots,a_{q-r})}\leq\cdots
\leq\frac{e_{r+1}(b_1,\ldots,b_q)}{e_{1}(a_1,\ldots,a_{q-r})}\leq
e_{r}(b_1,\ldots,b_q).
\end{equation}
Then the sequence of hypergeometric terms \emph{(}if $r=q$ the
numerator is $1$\emph{)},
$$
f_n=\frac{(a_1)_n\cdots(a_{q-r})_n}{(b_1)_n\cdots(b_q)_n},
$$
is log-concave, i.e. $f_{n-1}f_{n+1}\leq{f_n^2}$, $n=1,2,\ldots$
It is strictly log-concave unless $r=0$ and $a_i=b_i$,
$i=1,\ldots,q$.
\end{lemma}
The proof of this lemma for $r=0$ can be found in
\cite[Theorem~4.4]{HVV} and \cite[Lemma~2]{KS}. The latter
reference also explains how to extend the proof to general $r$
(see the last paragraph of \cite{KS}). This leads immediately to
the following statements.
\begin{theo}
Let $0\leq{p}\leq{q}$ be integers. Denote
$$
f(\nu,x):=\frac{1}{\Gamma(\nu)}{_pF_{q+1}}(a_{1},\ldots,a_{p};\nu,b_{1},\ldots,b_{q};x)
$$
and suppose that parameters $(a_1,\ldots,a_p)$, $(b_1,\ldots,b_q)$
satisfy \emph{(\ref{eq:symmetric-chain1})}. Then the function
$f(\nu,x)$ satisfies Theorem~\ref{th:gammadenom} and
Corollaries~\ref{cr:compl-mon}-\ref{cr:phi-below}.
\end{theo}
\begin{theo}
Let $0\leq{p}\leq{q+1}$ be integers. Denote
$$
g(\nu,x):=\frac{1}{\Gamma(\nu)}{_{p}F_{q+1}}(a_{1},\ldots,a_{p};\nu,b_{1},\ldots,b_{q};x)
$$
and suppose that parameters $(a_1,\ldots,a_p)$,
$(1,b_1,\ldots,b_q)$ satisfy \emph{(\ref{eq:symmetric-chain1})}.
Then the function $g(\nu,x)$ satisfies
Theorem~\ref{th:gammadenom1} and
Corollaries~\ref{cr:g-compl-mon}-\ref{cr:lambda-below}.
\end{theo}

\medskip

\textbf{Example~5}. Our last example is non-hypergeometric.
Consider the parameter derivative of the regularized Kummer
function:
$$
\frac{\partial}{\partial{a}}\frac{1}{\Gamma(b)}{_1F_1}(a;b;x)=
\sum\limits_{k=0}^{\infty}\frac{(\psi(a+k)-\psi(a))(a)_k}{\Gamma(b+k)k!}x^k
$$
(this function cannot be expressed as hypergeometric function of
one variable). If $a\geq{1}$ then the sequence
$$
h_k=\frac{(\psi(a+k)-\psi(a))(a)_k}{k!}
$$
is log-concave, since
\begin{multline*}
h_k^2-h_{k-1}h_{k+1}=\frac{(a)_{k-1}(a)_{k}}{(k-1)!k!}\times
\\
\left\{\frac{a+k-1}{k}(\psi(a+k)-\psi(a))^2-\frac{a+k}{k+1}(\psi(a+k-1)-\psi(a))(\psi(a+k+1)-\psi(a))\right\}>0.
\end{multline*}
The last inequality holds because $y\to{\psi(a+y)-\psi(a)}$ is
concave according to the Gauss formula~\cite[Theorem~1.6.1]{AAR}
$$
(\psi(a+y)-\psi(a))''_y=\psi''(a+y)=-\int\limits_{0}^{\infty}\frac{t^2e^{-t(a+y)}}{1-e^{-t}}dt<0
$$
and hence is log-concave while $(a+k-1)/k>(a+k)/(k+1)$ if
$a\geq{1}$. By Theorem~\ref{th:gammadenom1} this leads to
\begin{theo}
Suppose $a\geq{1}$.  Then the function
$$
h(\nu,x):=\frac{\partial}{\partial{a}}\frac{{_1F_1}(a;\nu;x)}{\Gamma(\nu)}
$$
satisfies Theorem~\ref{th:gammadenom1} and
Corollaries~\ref{cr:g-compl-mon}-\ref{cr:lambda-below}.
\end{theo}

\paragraph{5. Acknowledgements.} We acknowledge the
financial support of the Russian Basic Research Fund (grant
11-01-00038-a), Far Eastern Federal University and the Far Eastern
Branch of the Russian Academy of Sciences.

\end{document}